\documentclass[a4paper, 12pt]{article}
\usepackage{graphicx} 
\usepackage{setspace}
\usepackage{csquotes}
\usepackage[]{biblatex}
\usepackage{hyperref}

\usepackage[left=2cm,right=2cm,
    top=2cm,bottom=2cm,bindingoffset=0cm]{geometry}

\usepackage{amssymb,amsmath,amsthm} 

\usepackage{algorithm}
\usepackage{algpseudocode}

\newtheorem{statement}{Statement}
\newtheorem{remark}{Remark}

\usepackage[english]{babel}

\usepackage{mathrsfs}

\def\C{{\mathbb C}}
\def\R{{\mathbb R}}
\def\K{{\mathcal K}}
\def\L{{\mathcal L}}
\def\r{{\mathbf r}}
\def\v{{\mathbf v}}
\def\x{{\mathbf x}}
\def\y{{\mathbf y}}
\def\p{{\mathbf p}}
\def\u{{\mathbf u}}
\def\d{{\mathbf d}}
\def\e{{\mathbf e}}
\def\q{{\mathbf q}}
\def\0{{\boldsymbol \theta}}

\DeclareMathOperator{\LinSpan}{\mathscr{L}}
\DeclareMathOperator{\Colspan}{colspan}
\DeclareMathOperator{\im}{im}
\DeclareMathOperator{\Dim}{dim}

\newcommand{\norm}[1]{\left\lVert#1\right\rVert}
\renewcommand{\vec}[1]{\textbf{#1}}

\title{Generalized minimal residual method for systems with multiple right-hand sides}
\author{S.Sukmanyuk, D.Zheltkov, B.Valiakhmetov}
\date{}

\begin{document}
    
\onehalfspacing

\maketitle
\begin{abstract}

    We present a novel variant of the Generalized Minimal Residual method for linear systems with the same matrix and subsequent multiple right-hand sides.  The concept of the method is similar to that of the classical GMRES method: it minimizes the residual over the search space. However, the search space is no longer a Krylov subspace, since we do not restart the process after changing a right-hand side. Our way of the extension of the search space is arithmetically equivalent to the way of the Generalized Conjugate Residual method for multiple right-hand sides. However, our algorithm requires less memory, is less computationally complex, and maintains orthonormal direction vectors, which increases the robustness of the method.  Since there are no assumptions on the right-hand sides, the process can be viewed as an extension of the search space with an arbitrary vector. Therefore, our technique can be easily adapted for other GMRES variants such as flexible GMRES or GMRES with deflated restarting.
\end{abstract}

\section*{Keywords}
Generalized Minimal Residual method, Generalized Conjugate Residual method, Krylov subspace, multiple right-hand sides.

\section*{Introduction}

Most numerical methods for solving time-dependent partial differential equations lead to linear systems with the same non-singular matrix and multiple right-hand sides. All the right-hand sides are not necessarily available at once; e.g. another right-hand side may depend on the solution with the previous one. Such systems are usually solved with a direct LU-factorization. The matrix of the system has to be factorized once, which is the major advantage of this approach. However, direct solvers may require a lot of memory and be very computationally complex. Hence, iterative solvers are often more relevant for such problems. Moreover, they can take into account any information from operating with previous right-hand sides. For example, in the wave scattering problem~\cite{stavtsev2009application}, \cite{smith1989conjugate}, time-marching methods for PDEs~\cite{fischer1998projection}, or structural mechanics problems~\cite{farhat1994implicit}, the right-hand sides are not arbitrary. The most widely used iterative solvers are Krylov methods. We emphasize that in our case the right-hand sides are not available simultaneously, which prevents the use of Krylov block methods \cite{o1980block, robbe2006exact, baker2006improving, simoncini1995iterative}.

For the symmetric case, the extensions of the Lanczos and the Conjugate Gradient~(CG) method have been developed for multiple right-hand sides. As soon as another system is solved, a Galerkin projection onto the search space is performed to get approximations for the unsolved systems. If the desirable tolerance has not been achieved after the projection, two strategies for refinement are taken. The first is to start a new Lanczos or CG method with the obtained approximation as an initial guess. Methods with such a restarted Lanczos-Galerkin procedure, also known as seed methods, have been introduced in \cite{papadrakakis1990new}, \cite{van1987iterative}, \cite{smith1989conjugate} and analyzed in~\cite{chan1997analysis}. In reality, seed methods are effective if the right-hand sides are available simultaneously, so different strategies can be applied to choose a seed system. An alternative approach has been developed by Saad in~\cite{saad1987lanczos}, where the Lanczos run is continued instead of restarting the method for another right-hand side. The previous Lanczos process is continued naturally by the procedure introduced by Parlett~\cite{parlett1980new}. Compared to seed methods, this method does not require the right-hand sides to be simultaneously available, but tends to use more memory for the saved Lanczos vectors. 

For the non-symmetric case, our attention is restricted to the GMRES method developed by Saad and Schultz~\cite{saad1986gmres} for linear systems with a single right-hand side.  In fact, several extensions of the GMRES method have been developed to solve systems with multiple right-hand sides, most of which are close to seed methods. For example, in~\cite{simoncini1995iterative} and~\cite{nachtigal1992hybrid}  different techniques for the construction of better approximations to the unsolved systems are suggested. In~\cite{parks2006recycling} a way to recycle the search space is provided. However, all the methods lose information from operating with previous right-hand sides. 

An attempt is made in ~\cite{prasad1995gmres} to develop the GMRES method that extends the search space after a change in the right-hand side is introduced in~\cite{prasad1995gmres}. However, the relations presented contain inaccuracies; therefore, the method does not actually minimize the residual over a search space. Moreover, the method may break down for some systems. 

For all these reasons, our aim is to develop a method that allows a new right-hand side to be added to the GMRES process without restarting the method. The relations of the presented method are quite similar to the relations of the classical one: in every iteration $k$ an approximation to the solution of the $m$th system is chosen as $\x_k^{(m)} = \x_0^{(m)} + \y_k$, $\y_k \in \L_k$, where $\x_0^{(m)}$ is some initial question. The difference is that the search space $\L_k$ is no longer a Krylov space. However, our choice of search space is based on one of the representations of the Krylov subspace. We will show that, under some assumptions, the search space of our method is theoretically equivalent to the search space of the GCR method extended for multiple right-hand sides~\cite{lingen1999generalised}. We note that the classical GCR~\cite{eisenstat1983variational} method is theoretically equivalent to the classical GMRES method. However, the latter is less computationally complex, requires less memory, and is more robust. Our method preserves these advantages over the GCR method for multiple right-hand sides.

Actually, changing a right-hand side can be considered as an extension of the search space with an arbitrary direction vector. Therefore, our method can be easily adopted for GMRES-like methods with flexible directions such as flexible GMRES with variable preconditioner \cite{saad1993flexible, van1994gmresr, baglama1998adaptively} or deflated GMRES \cite{morgan2002gmres, giraud2010flexible, morgan1995restarted, morgan2000implicitly, chapman1997deflated}.

In Section 1 we recall the classical GMRES. Section 2 describes in detail the extension of the method for multiple right-hand sides. The main challenge in the construction of the method is the construction of the orthonormal basis of the search space, which is not a problem in the classical method. Therefore, a significant part of Section 2 is dedicated to it. We present the numerical results in Section 3, where we compare the new GMRES method to the GCR method for multiple right-hand sides and to the seed variant of GMRES for a wave-scattering problem.

\section{Classical GMRES method}
Let us recall the classical Generalized Minimal Residual (GMRES) methods. Our notation differs from the classical one. However, they are more convenient for further description of the method for multiple right-hand sides. 

The method iteratively solves a linear system $A \vec{x} = \vec{b}$ with a non-singular matrix $A \in \C^{N \times N}$ and a single right-hand side $\vec{b} \in \C^{N}$. An initial guess $\x_0 \in \C^N$ and a desired tolerance $\varepsilon \in \R$ are provided.  The initial residual is $\r_0 = \vec{b} - A\x_0$. On every iteration $k = 1, 2, \dots$ the method seeks an approximation to the solution in the affine Krylov subspace $\K_k(A, \r_0) = \LinSpan(\r_0, A\r_0, \dots, A^{k-1}\r_0)$:
\begin{equation}  \label{x_k}
    \x_k = \x_0 + \y_k, ~ \y_k \in \K_k(A, \r_0),
\end{equation} where $\y_k$ minimizes the Euclidean norm of the $k$th residual:
\begin{equation*}
\    \r_k = \vec{b} - A\x_k = \r_0 - A\y_k.
\end{equation*}
The norm of the residual is minimal when the orthogonal decomposition of $\r_0$ takes place:
\begin{equation}\label{r_0_orth_dec}
    \r_0 = \r_k +  A\y_k, ~\r_k \perp A\K_k(A, \r_0), ~A\y_k \in A\K_k(A, \r_0).
\end{equation}
For the representation of the solution as \eqref{x_k} the subspace $\K_k(A, \r_0)$ is called the \textit{search space}. Since the search space is a Krylov space, the Arnoldi process~\cite{arnoldi1951principle} is applied to construct its orthonormal basis. Let $Q_{k} = [\q_1, \dots, \q_{k + 1}]$ denote the orthonormal basis of $\K_{k + 1}(A, \r_0)$ such that $\q_1, \dots, \q_{l}$ forms the orthonormal basis of $\K_l(A, \r_0)$, $l=\overline{1, k + 1}$. The Arnoldi process provides the equation:
\begin{equation}\label{arnoldi}
    AQ_k\begin{bmatrix}
        I_k \\
        \0^T
    \end{bmatrix} = Q_{k}H_{k},
\end{equation}
where $H_{k} \in \C^{(k+1) \times k}$ such that $\{H_k\}_{ij} = 0$ for $i  > j + 1$. Here $Q_k\begin{bmatrix}
        I_k \\
        O
    \end{bmatrix} = [\q_1, \dots, \q_k]$  is actually the orthonormal basis of $\K_k(A, \r_0)$.

As soon as we construct $Q_{k}$, we represent $\y_k$  as $\y_k = Q_{k}\begin{bmatrix}
        I_k \\
        \0^T
    \end{bmatrix} \v_k$, where $\v_k \in \C^k$. In addition, from \eqref{arnoldi} it is clear that $\im(Q_kH_k) = A\K_k$. With QR factorization of $H_k$: $H_k = G_k\begin{bmatrix}
        R_k \\
        \0^T
    \end{bmatrix}$, where $G_k \in \C^{(k + 1) \times (k + 1)}$ is a unitary matrix, $R_k \in \C^{k \times k}$ is an upper-triangular one, we get  $U_k = Q_kG_k \begin{bmatrix}
        I_k \\
        \0
    \end{bmatrix}$ is the orthonormal basis of $A\K_k(A, \r_0)$. Having the representation of $U_k$ we easily transform \eqref{r_0_orth_dec} to:
\begin{align*}
    \v_k &= R_k^{-1}U_k^*\r_0 = R_k^{-1}\begin{bmatrix}
        I_k & \0_k \end{bmatrix}G_{k+1}^*Q_{k+1}^*\r_0\\
    \r_k &= (I - U_kU_k^*)\r_0 =  \r_0 - Q_{k+1}G_{k+1} \begin{bmatrix}
    I_k  & \0_k\\
    \0_k^T & 0
\end{bmatrix} G_{k+1}^*Q_{k+1}^*\r_0.
\end{align*}
Since $\r_0 \in \K_{k+1} = \Colspan(Q_{k+1})$, it can be represented as $\r_0 = Q_{k+1}Q_{k+1}^*\r_0$, resulting in  
\begin{equation}\label{residual}
    \r_k =  Q_{k+1}G_{k+1}\begin{bmatrix}
        O_k  & \0_k\\
        \0_k^T & 1
    \end{bmatrix} G_{k+1}^*Q_{k+1}^*\r_0.
\end{equation}
Moreover, there is no need to explicitly calculate the residual to compute its norm. Indeed, since \eqref{residual} holds and the Euclidean norm is unitary-invariant, then
\begin{equation}
\label{eq:r_k_norm_estimate}
    \norm{\r_k} = | \vec{e}_{k+1}^TG_{k+1}^*Q_{k+1}^*\r_0|,
\end{equation}
where $\vec{e}_{k+1} = [0 ~  \dots ~ 0 ~ 1] \in \C^{k+1}$.
So in every iteration the residual norm is estimated as \eqref{eq:r_k_norm_estimate}. The solution $\v_k$ is computed when the desired tolerance is achieved: $\norm{\r_k} \le \varepsilon \norm{\vec{b}}$.

The final procedure of the GMRES method is presented in Algorithm~\ref{alg:GMRES}.

To describe how to extend the method for multiple right-hand sides, we emphasize several properties of the classical method.
\begin{enumerate}
    \item An approximation to the solution is actually sought as
    $$
        \x_k = \x_0 + P_k \v_k,~ \v_k \in \C^k
    $$
    where $P_k = Q_k \begin{bmatrix}
        I_k \\
        \0^T
    \end{bmatrix}$. With this formulation, the columns of $P_k$ are called \textit{direction vectors}. We note that the direction vectors are orthonormal.

    \item If we denote the search space as $\L_k$: $\L_k = \K_{k}(A, \r_0)$ we can note that $\L_k, A\L_k \in \L_{k+1}$ therefore $Q_{k}$ is a common basis for $\L_k$ and $A\L_k$. Therefore the corresponding bases $P_k$ and $U_k$ are effectively represented within $Q_{k}.$

    \item  Let us denote $\u_1, \dots, \u_k$ the orthonormal basis of $A\L_k$ such that $\u_1, \dots, \u_{s}$ is the orthonormal basis of $A\L_s$ for $s = \overline{1, k}$. Since $\L_k = \K_k(A, \r_0)$, it can be represented as
    \begin{equation}\label{GMRES_ss}
        \L_k = \LinSpan(\r_0) + A\L_{k-1} = \LinSpan(\r_0, \u_1, \dots, \u_{k-1}) = \L_{k-1} + \LinSpan(\u_{k-1}).
    \end{equation}
\end{enumerate}

\begin{algorithm}[ht]
 \begin{algorithmic}
 \Require $A$, $\vec{b}$, $\x_0$, $\varepsilon$
\State  $k = 0$
\State $\r_0 = \vec{b} - A\x_0$
\State $\q_1 = \r_0 / \norm{\r_0 }$
\State $Q = [\q_1]$
\State $H = [~]$
\State $G = [1]$
\While {$ \norm{\r_0}|\vec{e}_{k+1}^TG^*\vec{e}_{k+1}| > \varepsilon\norm{\vec{b}} $}
\State $k = k + 1$
\State $\q_{k+1} = A\q_{k}$
\State Arnoldi($\q_{k+1}$, $Q$, $H$)
\State  $G, R$ = qr($H$)
\EndWhile
\State $\vec{v} = 
     R^{-1}\begin{bmatrix}
        I_k  & \0\\
    \end{bmatrix} G^*Q^*\r_0$
\State $\x = \x_0 + Q \begin{bmatrix}
   \vec{v}\\
   0
\end{bmatrix}$
\end{algorithmic}
\caption{GMRES method}
\label{alg:GMRES}
\end{algorithm}

\begin{algorithm}[ht]
\caption{Arnoldi}
\label{alg:arnoldi}
 \begin{algorithmic}
 \Require $\q$, $Q$, $H$
\State $\vec{h} = Q^*\vec{q}$
\State $\vec{q} = \q - QQ^*\vec{h}$
\State $H = \begin{bmatrix}
    H & \vec{h}\\
    \0^T & \norm{\q}
\end{bmatrix}$
\State $\q = \q / \norm{\q}$
\State $Q = \begin{bmatrix}
    Q & \q
\end{bmatrix}$
\end{algorithmic}
\end{algorithm}

\section{GMRES method for multiple right-hand sides}

\subsection{General idea}

Our aim is to extend the GMRES method for linear systems with a non-singular  matrix $A \in \C^{N \times N}$  and multiple right-hand sides $\vec{b}^{(s)} \in \C^{N}$, $s=\overline{1, M}$:
\begin{equation*}
    A\x^{(s)} = \vec{b}^{(s)},~ s \in \overline{1, M}.
\end{equation*}
The initial guess $\x^{(1)}_0, \dots, \x^{(M)}_0$ and the required tolerance $\varepsilon^{(1)}, \dots, \varepsilon^{(M)}$ are provided.  The initial residuals are $\r^{(s)}_0 = \vec{b}^{(s)} - A\x^{(s)}_0$, $s \in \overline{1, M}$. The right-hand sides are not assumed to be available simultaneously. For example, the next right-hand side may depend on the solution of the previous system. Therefore, systems are solved sequentially.

The general concept of the method is the same as in the classical one: solving the minimization problem in the search space. The only difference is that in the $k$th iteration a certain right-hand side is being processed. Let us denote $k$ as the current iteration and $m$ as a right-hand side number that is being processed in the $k$th iteration (we imply that the $1$st, $2$nd, \dots, $(m - 1)$th right-hand sides have already converged in previous iterations).

By analogy to the classical method on the $k$th iteration, the search space $\L_k$ is chosen (actually $\L_1 \subset \L_2 \subset \dots \subset \L_k$ and $\Dim \L_j = j $). The dimension of the search space does not depend on the right-hand side number.

An approximation to the solution for the $m$th right-hand side on the $k$th iteration is 
\begin{equation}\label{x_k_multiple}
    \x_k^{(m)} = \x_0^{(m)} + \y_k^{(m)},
\end{equation}
where $\y_k^{(m)} \in \L_k$ minimizes the norm of the residual:
\begin{equation}\label{r_k_multiple}
    \r_k^{(m)} = \vec{b}^{(m)} - A\x_k^{(m)} = \r_0^{(m)} - A\y_k^{(m)}.
\end{equation}

The residual is minimal when the orthogonal decomposition takes place:
\begin{equation}\label{r_k_orth_dec_multiple}
    \r_0^{(m)} =  \r_k^{(m)}  + A\y_k^{(m)}, ~ \r_k^{(m)} \perp A\L_k, ~A\y_k^{(m)} \in A\L_k.
\end{equation}

If the desired tolerance is achieved on the $k$th iteration: $\norm{\r_k^{(m)}} \le \varepsilon^{(m)} \norm{\vec{b}^{(m)}}$, the method changes the right-hand side. For the $(m + 1)$th right-hand side we can construct an approximate solution using $\L_k$ the same as \eqref{x_k_multiple}, \eqref{r_k_multiple}: $\r_k^{(m+ 1)}$ is orthogonal to the $A\L_k$ part of $\r_0^{(m + 1)}$. Without loss of generality, we assume that the desired tolerance is never achieved in $\r_k^{(m + 1)}$: $\norm{\r_k^{(m + 1)}} > \varepsilon^{(m + 1)} \norm{\vec{b}^{(m + 1)}}$. So in the next iteration, the right-hand side is $m + 1$.

For every right-hand side $s$ we denote $n_s$ as the iteration number on which the $s$th right-hand side converged. We set $n_0 = 0$.

To construct the method for multiple right-hand sides, we need to specify the search space $\L_k$. It is clear that it cannot be a Krylov subspace any longer since we operate on several right-hand sides.

\subsection{The search space}
Our choice of search space is based on the representation \eqref{GMRES_ss} of the search space for a single right-hand side. To describe it, we recall that $k$ denotes the iteration number, $m$ denotes the right-hand side that is being processed on the $k$th iteration. We will describe the construction of the search space $\L_k$ by induction. We set $\L_0 = \LinSpan(\0)$.

Assume that $\L_{k-1}$ and the orthonormal basis $\u_1, \dots, \u_{k-1}$ of $A\L_k$ have been defined. Two cases take place.
\begin{enumerate}
    \item  The $(m-1)$th right-hand side converged on the previous iteration: $n_{m-1} = k - 1$. The residual for the $m$th right-hand side $\r_{n_{m-1}}^{(m)}$ is orthogonal to the $A\L_{n_{m-1}}$ part of $\r_0^{(m)}$.  In this case $\r_{n_{m-1}}^{(m)}$ extends the search space:
    $$
        \L_{k} = \L_{k-1} + \LinSpan(\r_{n_{m-1}}^{(m)}).
    $$
    It has been assumed that $\norm{\r_{n_{m-1}}^{(m)}} > \varepsilon^{(m)} \norm{\vec{b}^{(m)}}$. 
    \item The previous iteration was performed on the $m$th right-hand side.  In this case $\u_{k-1}$ extends the search space:
    $$
        \L_{k} = \L_{k-1} + \mathscr{L}(\u_{k - 1}).
    $$
\end{enumerate}
For both cases $\u_{k}$ is defined as $\u_1, \dots, \u_{k-1}, \u_k$ form the orthonormal basis of $A\L_k$.

According to the process of the construction, the search space $\L_k$ is represented as:
\begin{align}\label{search_space}
    \L_{k} &= \LinSpan(\r_0^{(1)}, \u_1, \dots, \u_{n_1 - 1}) + \\ \nonumber
    &+ \LinSpan(\r_{n_1}^{(2)}, \u_{n_1 + 1}, \dots, \u_{n_2 - 1}) +  \\\nonumber
    &+ \dots + \\ \nonumber
    &+ \LinSpan(\r_{n_{m-1}}^{(m)}, \u_{n_{m-1} + 1}, \dots, \u_{k - 1}),
\end{align}
where $A\L_j = \LinSpan(\u_1, \dots, \u_j)$, $j = \overline{1,k}$.

\begin{remark}

     Since $\r_j^{(s)}$ is orthogonal to the $A\L_j$ part of $\r_0^{(s)}$ for $s = \overline{1, m}$ and corresponding $j$, the search space \eqref{search_space} is equal to 
      \begin{align*}\L_{k} &= \LinSpan(\r_0^{(1)}, \r^{(1)}_1, \dots, \r^{(1)}_{n_1 - 1}) + \\ \nonumber
    &+ \LinSpan(\r_{n_1}^{(2)}, \r^{(2)}_{n_1 + 1}, \dots, \r^{(2)}_{n_2 - 1}) +  \\\nonumber
    &+ \dots + \\ \nonumber
    &+ \LinSpan(\r^{(m)}_{n_{m-1}}, \r^{(m)}_{n_{m-1} + 1}, \dots, \r^{(m)}_{k - 1})
    \end{align*}
    if $(\r_{j}^{(s)}, \u_{j}) \ne 0$ for iterations $j$ performed with the $s$th right-hand side.   This subspace is chosen for the extended GCR method for multiple right-hand sides \cite{lingen1999generalised}. The classical GCR method \cite{eisenstat1983variational} is theoretically equivalent to the classical GMRES method, which in turn is more robust, requires less memory and has less numerical complexity. However, the GCR method can be easily extended for multiple right-hand sides.
\end{remark}

In the next Section, we show that for our choice of the search space \eqref{search_space}, $\Dim \L_k = k$, thus the method provides the solution in at most $N$ iterations.

\subsection{Analysis of the search space}

For $\L_k$ constructed as \eqref{search_space}, several statements hold. Let us recall that $k$ denotes the iteration, $m$ denotes the right-hand side that is being processed in the $k$th iteration (meaning that the right-hand sides $1, \dots, (m-1)$ have converged in previous iterations). And $n_s$ denotes the iteration when the $s$th right-hand side converged, $s = \overline{1, m - 1}$. Vectors $\u_1, \dots \u_j$ form the orthonormal basis of $A\L_j$ for $j = \overline{1, k}$.
\begin{statement}\label{st_1}
    For every right-hand side $ 1 \le s \le m$:
    \begin{equation}
        \r^{(s)}_{j} = \r_{j - 1}^{(s)} - (\r_{j-1}^{(s)}, \u_j)\u_j,
    \end{equation}
    where $j = \overline{n_{(s-1)}, n_s}$ for $s = \overline{1, m - 1}$ and $j = \overline{n_{(m-1)}, k}$ for $s = m$.
\end{statement}
\begin{proof}
    It is clear since $\r^{(s)}_{j}$ is orthogonal to the $A\L_j = \LinSpan(\u_1, \dots, \u_j)$ part of $\r_0^{(s)}$.
\end{proof}

\begin{statement}\label{st_2}
For every right-hand side $ 1 \le s \le m - 1$:  
\begin{equation}
    (\r_{n_{s}-1}^{(s)}, \u_{n_s})   \ne 0,
\end{equation}
where $n_{s}$ denotes the iteration number when the $s$th right-hand side converged:$\norm{\r_{n_{s}}^{(s)}} \le \varepsilon \norm{\vec{b}^{(s)}}$.
\end{statement}
\begin{proof}

    \item  The $s$th system converged in the $n_s$th iteration : the residual norm has decreased from the previous iteration. This means that $\norm{\r_{n_{s}}^{(s)}} \le \varepsilon \norm{\vec{b}^{(s)}} < \norm{\r_{n_{s}-1}^{(s)}}$. Using statement \ref{st_1}:
    $$
        \r^{(s)}_{n_s} = \r^{(s)}_{n_s-1} - (\r^{(s)}_{n_s - 1}, \u_{n_s}) \u_{n_s},
    $$
    we get $(\r^{(s)}_{n_s - 1}, \u_{n_s}) \ne 0$. But $  (\r^{(s)}_{n_{s}-1}, \u_{n_s}) =  (\r^{(s)}_{n_{(s-1)}}, \u_{n_s}) = (\r^{(s)}_{0}, \u_{n_s})$, since $\r^{(s)}_{n_{(s-1)}}$ and $\r^{(s)}_{n_{s}-1}$ are orthogonal  projections of  $\r^{(s)}_{0}$  to $(A\L_{n_{(s-1)}})^{\perp}$ and $(A\L_{n_{s} - 1})^{\perp}$ correspondingly.
\end{proof}
\begin{statement}\label{st_3}
In the $k$th iteration $\Dim \L_k = k$.
\end{statement}
\begin{proof}
For $k = 1$ $\L_k = \LinSpan(\r_0^{(1)})$. The initial residual $\r_0^{(1)}$ is considered to be non-zero, so $\Dim \L_1 = 1$.

For $k > 1$ assume that $\Dim \L_{k-1} = k-1$. There are two cases.
    \begin{enumerate}
    \item $k$ is the first iteration with the current system: $n_{m-1} = k - 1$ and $\L_{k} = \L_{n_{m-1}} + \LinSpan(\r^{(m)}_{n_{m-1}})$.
    
    If $\Dim \L_k \ne k$, then  \[\r^{(m)}_{n_{m-1}} \in \L_{n_{m-1}} = \LinSpan(\r_0^{(1)}, \u_1, \dots, \u_{n_1 - 1}) + \dots + \LinSpan(\r_{n_{m-2}}^{(m-1)}, \u_1, \dots, \u_{n_{m-1} - 1}).\]
    So the residual $\r^{(m)}_{n_{m-1}}$ can be represented as
        \begin{align*}
             \r^{(m)}_{n_{m-1}} &= \sum_{s=1}^{m-1} \alpha_s \r_{n_{s-1}}^{(s)} +                     \sum_{j=1}^{n_{m-1}} \beta_j \u_j,
        \end{align*}
        where $\beta_{n_1} = \beta_{n_2} = \dots = \beta_{n_{m-1}} = 0$.

        Since statement~\ref{st_1} holds, the dot product of  $\r^{(m)}_{n_{m-1}}$ and $\u_{n_s}$ for $s=\overline{1, m - 1}$ transforms to
    $$
       0 = (\r^{(m)}_{n_{m-1}}, \u_{n_s}) = \alpha_1 (\r_{0}^{(1)}, \u_{n_s}) + \alpha_2 (\r_{n_1}^{(2)}, \u_{n_s}) + \dots + \alpha_s (\r_{n_{s-1}}^{(s)}, \u_{n_{s}}) .
    $$
    Since this holds for $s=\overline{1, m - 1}$ and from statement~\ref{st_2}  we have $(\r_{n_{s-1}}^{(s)}, \u_{n_s}) \ne 0$, we sequentially obtain  $\alpha_1 = \alpha_2 = \dots = \alpha_{m - 1} = 0$.  Having this, it is clear that $\beta_1 = \dots = \beta_{n_{m-1}} = 0$.  Finally, $\r^{(m)}_{n_{m-1}} = \0$ that was excluded. So $\r^{(m)}_{n_{m-1}} \not \in \L_{k-1}$  meaning that $\Dim \L_{k} = \Dim \L_{k - 1} + 1 = k$.
    
    \item If $k$ is not the first iteration with the current system, then $\L_{k} = \L_{k-1} + \LinSpan(\u_{k-1})$. If $\Dim \L_{k} \ne k$, then
    \begin{align*}
            \u_{k-1} &\in \L_{k-1} = \LinSpan(\r^{(1)}_0, \dots,\r_{n_{m-1}}^{(m)}) + \LinSpan( \u_1, \dots, \u_{n_1 - 1}, \u_{n_1 + 1}, \dots, \u_{k-2}), \\
             \u_{k-1} &= \sum_{s=1}^{m} \alpha_s \r_{n_{s-1}}^{(s)} +                     \sum_{j=1}^{k-2} \beta_j \u_j,
    \end{align*}
    where $\beta_{n_1} = \beta_{n_2} = \dots = \beta_{n_{m-1}} = 0$.

    The same way as in the previous case, we can obtain that $\alpha_1 = \dots = \alpha_{m-1} = 0$ and $\beta_1 = \dots = \beta_{n_{m-1}} = 0$ therefore 
    \begin{align*}
             \u_{k-1} &= \alpha_m \r_{n_{m-1}}^{(m)} +                     \sum_{j=n_{m-1} + 1}^{k-2} \beta_j \u_j.
    \end{align*}
    It means that $\alpha_m\r_{n_{m-1}}^{(m)} \in \mathscr{L}(\u_1, \dots, \u_{k-1}) = A\L_{k-1}$.  But $\r_{k-1}^{(m)}$ is a projection of $\r_{n_{m-1}}^{(m)}$ to $( A\L_{k-1})^{\perp}$. On the previous iteration the required tolerance was not achieved, therefore $\r_{n_{m-1}}^{(m)} \not \in A\L_{k-1}$ and $\alpha_m = 0$. Since $\alpha_m = 0$, it is obvious that $\beta_j = 0$, $j = \overline{n_{m-1} + 1, \dots, k-2}$. Hence $\u_{k-1} = \0$ that conflicts with $\Dim \L_{k-1} = k - 1$.
    \end{enumerate}
\end{proof}

The search spaces are nested: $\L_1 \subset \L_2 \subset \dots \subset \L_{k-1} \subset \L_{k}$ and $\Dim \L_j = j$ for $j = \overline{1, k}$.

\subsection{Common basis}
From \eqref{search_space} we note that
\begin{equation}\label{common_ss_1}
    \L_k + A\L_k = \LinSpan(\r_0^{(1)}, \dots, \r^{(m)}_{n_{m-1}}) + \LinSpan(\u_1, \dots, \u_k).
\end{equation}
Thus $\Dim(\L_k + A\L_k) \le k + m$, therefore we can construct the basis of $\L_k + A\L_k$ and represent the basis $U_k$ and the direction vectors $P_k = [\p_1, \dots, \p_k]$ within it.

We can simplify the representation \eqref{common_ss_1} recalling that $\r^{(s)}_{n_{s-1}}$ is a projection of $\r_0^{(s)}$ to $(A\L_{n_{s-1}})^{\perp} = \LinSpan(\u_1, \dots, \u_{n_{s-1}})^{\perp}$, $n_{s - 1} < k$, $s=\overline{1,m}$. We obtain
\begin{equation}\label{L_comm}
    \L_k + A\L_k = \LinSpan(\r_0^{(1)}, \dots, \r^{(m)}_{0}) + \LinSpan(\u_1, \dots, \u_k) = \LinSpan(\r_0^{(1)}, \dots, \r^{(m)}_{0}) + \LinSpan(A\p_1, \dots, A\p_k).
\end{equation}
Similarly we note 
\begin{equation}\label{L_comm_hat}
    \L_k + A\L_{k - 1} = \LinSpan(\r_0^{(1)}, \dots, \r^{(m)}_{0}) + \LinSpan(\u_1, \dots, \u_{k - 1}) = \LinSpan(\r_0^{(1)}, \dots, \r^{(m)}_{0}) + \LinSpan(A\p_1, \dots, A\p_{k - 1}).
\end{equation}
Let us denote
\begin{itemize}
    \item $t_k = \Dim (\L_k + A\L_{k})$, $s_k = \Dim (\L_k + A\L_{k - 1})$;
    \item $Q_{k} \in \C^{N \times t_k}$ the orthonormal basis of $\L_k + A\L_k$;
    \item $\hat{Q}_{k} \in \C^{N \times s_k}$ the orthonormal basis of $\L_k + A\L_{k - 1}$. From \eqref{L_comm} and \eqref{L_comm_hat} we see that if $k > n_{m - 1} + 1$ then $\hat{Q}_{k} = Q_{k-1}$ and $s_k = t_{k - 1}$.
\end{itemize}

\begin{enumerate}
    \item If $k = n_{m - 1} + 1$ then $Q_{k - 1}$ is appended with $\r_0^{(m)}$ and $A\p_k$ to get $Q_k$ and   $Q_{k - 1}$ is  appended with $\r_0^{(m)}$ to get $\hat{Q}_k$.

    \item If $k = n_{m - 1}$ then $Q_k$ is $Q_{k - 1}$ appended with $A\p_k$, $\hat{Q}_k$ is $Q_{k - 1}$. And for every $k$ we have $\hat{Q}_{k} = Q_{k}\begin{bmatrix}
        I_{s_k} \\
        O
    \end{bmatrix}$
\end{enumerate}
For construction of $Q_{k}$, we need the direction vector $\p_k$ to be available in the $k$th iteration. At the same time, to construct $\hat{Q}_k$, we do not need $\p_k$ to be available. However, $\p_k \in \L_k \subset \hat{Q}_k$ means that on every iteration $k$ the direction vector $\p_k$ can be computed through $\hat{Q}_k$:
$$
     \p_k = \hat{Q}_k \vec{c}_k,~ \vec{c}_k \in \C^{s_k}.
$$
As soon as $\p_k$ is available, $Q_{k}$ can be constructed. As soon as we have $Q_{k}$, all direction vectors can be represented through it:
$$
    P_k = Q_k C_k, ~C_k \in \C^{t_k \times k}.
$$
If we do not store $P_k$, we need to store the matrix $C_k$ to compute the solution $\x_k^{(m)} = \x_0^{(m)} + AP_k\v_k^{(m)}$ when $k = n_m$. It is also clear how $\vec{c}_k$ corresponds to the last column of $C_k$. The next section describes how $\vec{c}_k$ and $C_k$ are constructed.

Since $\Colspan(AP_j) = A\L_j = \Colspan(Q_{j})$, $j = \overline{1, k}$, we can get the extended Arnoldi equation:
\begin{eqnarray} \label{Arnoldi_advanced}
    A P_k = Q_{k} H_{k},
\end{eqnarray}
$ H_{k} \in C^{t_k \times k}$ and $[H_{k}]_{i, j} = 0$ for $i > t_j$.

It is clear from \eqref{Arnoldi_advanced} that $Q_{k}H_{k}$ is the basis of $A\L_k$.  A QR factorization of $H_{k}$ is performed to construct its orthonormal basis: 

$$H_{k} = G_{k}\begin{bmatrix} R_k \\ O \end{bmatrix},$$ where $G_{k} \in \C^{t_k \times t_k}$ is a unitary matrix, $R_k \in \C^{k\times k}$ is an upper triangular one. 

We note that $U_k = Q_{k}G_{k} \begin{bmatrix}
    I_k \\
    O
\end{bmatrix}$ is the orthonormal basis of $A\L_{k}$ such that $U_j = Q_{j}G_{j} \begin{bmatrix}
    I_j \\
    O\end{bmatrix}$ forms an orthonormal basis of $A\L_{j}$, $j=\overline{1,k}$. The calculation of $\r_k$ and $\v_k$ from \eqref{r_k_orth_dec_multiple} is reduced to:

$$
   \v_k^{(m)} =  R_k^{-1}U_k^*\r_0= R_k^{-1}\begin{bmatrix}
       I_{k} & O
   \end{bmatrix}G_{k}^*Q_{k}^*\r_0^{(m)},
$$

$$\r_{k}^{(m)} = (I - U_kU_k^*)\r_{0}^{(m)} = \Big(I - Q_{k}G_{k}\begin{bmatrix}
    I_k & O \\
    O & O_{t_k - k}
\end{bmatrix} G_{k}^*Q_{k}^*\Big)\r_{0}^{(m)}.$$
As $\r_{0}^{(m)} \in \L_{k}+A\L_k = \Colspan(Q_{k})$, then $\r_{0}^{(m)} = Q_{k}Q_{k}^*\r_{0}^{(m)}$. $G_{k}$ is a unitary matrix:
\begin{equation}\label{r_k_m}
    \r_{k}^{(m)} = Q_{k}G_{k}\begin{bmatrix}
    O_k & O \\
    O & I_{t_k - k}
\end{bmatrix} G_{k}^*Q_{k}^*\r_{0}^{(m)}.
\end{equation}
There is no need to compute the residual to estimate its norm:
$$
    \norm{\r_k^{(m)}} = \norm{\begin{bmatrix}
    O & I_{t_k - k}
\end{bmatrix} G_{k}^*Q_{k}^*\r_{0}^{(m)}}.
$$

It remains to determine how to construct the direction vector $\p_k$ without explicit storage of the entire basis $P_k$.

\subsection{Direction vectors}

In the current section, we describe how to construct the next direction vector $\p_{k+1}$ on the $(k+1)$th iteration. From the $k$th iteration, we have the matrix $Q_{k}$. However, in the $(k + 1)$th iteration, we can construct $\hat{Q}_{k+1}$ before $\p_{k+1}$. Recalling that $\p_{k + 1} \in \L_{k + 1} \subset \Colspan(\hat{Q}_{k + 1})$, we can seek the representation of $\p_{k + 1}$ as
$$
    \p_{k+1} = \hat{Q}_{k+1}\vec{c}_{k+1}, ~\vec{c}_{k+1} \in \C^{s_k}.
$$
Let us introduce the following auxiliary notation:
\begin{enumerate}
    \item $\widetilde{\p}_{k+1} = \u_k$ and  $\hat{G}_{k+1} = G_k \in \C^{s_{k+1} \times s_{k+1}}$ for  $k \ne n_{m}$,
    \item $\widetilde{\p}_{k+1} = \r_{n_{m}}^{(m)}$ and  $\hat{G}_{k+1} = \begin{bmatrix}
        G_k & \0 \\
        \0^T & 1
    \end{bmatrix} \in \C^{s_{k+1} \times s_{k+1}}$ if $k = n_{m}$,

    \item  $\Hat{E} = \begin{bmatrix}
        \e_{n_1} & \e_{n_2} & \dots & \e_{n_{m-1}} & \e_k & \e_{k+1} & \dots & \e_{t_k}
    \end{bmatrix} \in \C^{t_k \times (m + t_k - k)}$.
\end{enumerate}

To construct the direction vector $\p_{k+1}$, we recall that the direction vectors are orthonormal. Therefore, we need to orthogonalize $\widetilde{\p}_{k+1}$ to $\L_k$:
\begin{equation*}
    \L_k = \LinSpan(\u_1, \dots, \u_{n_1 - 1}, \u_{n_1 + 1}, \dots, \u_{n_{m-1} + 1}, \dots, \u_{k-1}) + \LinSpan(\r^{(1)}_0, \dots, \r^{(m)}_{n_{(m-1)}}).
\end{equation*}
In the previous section, we obtained $ [\u_1, \dots, \u_k] = U_k = Q_kG_k \begin{bmatrix}
    I_k \\
    O
\end{bmatrix}$. We can transform it to $U_k = \hat{Q}_{k+1}\hat{G}_{k+1} \begin{bmatrix}
    I_k \\
    O
\end{bmatrix}$. As a result,
\begin{equation}\label{L_k_QG}
    \L_k = \hat{Q}_{k+1}\hat{G}_{k+1}(\LinSpan(\e_1, \dots, \e_{n_1 - 1}, \e_{n_1 + 1}, \dots, \e_{n_{m-1} + 1}, \dots, \e_{k-1}) + \LinSpan(\d^{(1)}, \dots, \d^{(m)})),
\end{equation}
 where $\e_j \in \C^{s_{k+1}}$ is the $j$th element of the natural basis,   $\d^{(s)} = \hat{G}_{k+1}^*\hat{Q}_{k+1}^*\r^{(s)}_{n_{s-1}} \in \C^{s_{k+1}}$, $s=\overline{1, m}$. Let $D$ denote $D = [\d^{(1)}, \dots, \d^{(m)}]$.

 Since $\widetilde{\p}_{k+1} \in \L_{k+1} \subset \Colspan(\hat{Q}_{k+1}\hat{G}_{k+1})$, and $\hat{Q}_{k+1}\hat{G}_{k+1}$ is a unitary matrix, we can orthogonalize $\widetilde{\p}_{k+1}$ to $\L_k$ within the basis $\hat{Q}_{k+1}\hat{G}_{k+1}$. Considering \eqref{L_k_QG}, it is clear that we can perform the orthogonalization on the basis $\hat{Q}_{k+1}\hat{G}_{k+1}\hat{E}$. So, if $D_Q$ denotes an orthonormal basis of $\hat{E}^*D$, then the next direction vector is
 \begin{align}
     \hat{\p}_{k+1} &= \hat{Q}_{k+1}\hat{G}_{k+1}\hat{E}(I - D_QD_Q^*)\hat{E}^*\hat{G}_{k+1}^* \hat{Q}_{k+1}^*\widetilde{\p}_{k+1},\\
     \p_{k+1} &= \hat{\p}_{k+1} / \norm{\hat{\p}_{k+1}}
 \end{align}

The basis $D_Q$ can be obtained by QR factorization of $\hat{E}D \in \C^{(m + t_k - k) \times m}$. 
\begin{remark}
    We denote $\d^{(s)} = \hat{G}_{k+1}^*\hat{Q}_{k+1}^*\r^{(s)}_{n_{s-1}}$, $s=\overline{1, m}$. Recalling that $\r_{n_{s - 1}}^{(s)} = (I - U_{n_s}U_{n_s}^*)\r_{0}^{(s)}$, $U_{n_s} = \hat{Q}_{k+1}\hat{G}_{k+1}\begin{bmatrix}
    I_{n_s} \\
    O
\end{bmatrix}$ and $\r_0^{(s)} \in \Colspan(\hat{Q}_{k+1}\hat{G}_{k+1})$ we obtain
     $ \d^{(s)} = \begin{bmatrix}
     O_{n_{s-1}} \\
     I_{s_{k+1} - n_{s-1}}
 \end{bmatrix}\hat{G}_{k+1}^*\hat{Q}_{k+1}^*\r^{(s)}_{0}$. For the same reason,  $ \d^{(m)} = \begin{bmatrix}
     O_{k} \\
     I_{s_{k+1} - k}
 \end{bmatrix}\hat{G}_{k+1}^*\hat{Q}_{k+1}^*\r^{(m)}_{0}$. And considering that $\hat{E}^*D$ consists of corresponding components of $\d^{(1)}, \dots, \d^{(m)}$ we get that the matrix $\Hat{E}^*D$ is actually lower-triangular:
$$
\Hat{E}^*D = \begin{bmatrix}
    d^{(1)}_{n_1} & 0                   & \dots  &  0                        & 0 \\
    d^{(1)}_{n_2}       & d^{(2)}_{n_2} & \dots  &  0                        & 0 \\
    \vdots              & \vdots              & \vdots & \vdots                    & \vdots \\
    d^{(1)}_{n_{m-1}}   & d^{(2)}_{n_{m-1}}   & \dots  & d^{(m-1)}_{n_{m-1}} &  0 \\
    d^{(1)}_{k}         & d^{(2)}_{k}         & \dots  & d^{(m-1)}_k               & d^{(m)}_k \\
    d^{(1)}_{k+1}       & d^{(2)}_{k+1}       & \dots  & d^{(m-1)}_{k+1}           & d^{(m)}_{k+1} \\
    \vdots              & \vdots              & \vdots & \vdots                    &  \vdots \\
 d^{(1)}_{s_{k+1}}  & d^{(2)}_{s_{k+1}}       & \dots  & d^{(m-1)}_{s_{k+1}}           & d^{(m)}_{s_{k+1}} \\
\end{bmatrix}.
$$

Along with it, the last column $\Hat{E}^*\vec{d}^{(m)} = \Hat{E}^*\hat{G}_{k+1}^*\hat{Q}_{k+1}^*\r^{(m)}_{n_{m-1}}$ is not null as in the other case  the exact solution for the $m$th right-hand side was reached on the previous $(k-1)$th iteration. Recalling \eqref{st_1}, we obtain $d_{n_s}^{(s)} = (\r_{n_{s-1}}^{(s)}, \u_{n_s}) \ne 0$, $s = \overline{1, m-1}$. As a result, the columns of $\Hat{E}^*D$ are linearly independent. It again proves that $\Dim \L_k = k$.
\end{remark}

\begin{remark}
\begin{enumerate}
    \item If $k \ne n_{m}$, then $\widetilde{\p}_{k+1} = \u_k$ and $\hat{E}^*\hat{G}_{k+1}^*\hat{Q}_{k+1}^*\widetilde{\p}_{k+1} = \e_{k} \in \C^{m + t_k - k}$.

    \item If $k = n_{m}$ then $\widetilde{\p}_{k+1} = \r_{n_m}^{(m+1)}$ and $\hat{E}^*\hat{G}_{k+1}^*\hat{Q}_{k+1}^*\widetilde{\p}_{k+1} = \begin{bmatrix}
     O_{m} \\
     I_{t_k - k}
 \end{bmatrix}\hat{E}^*\hat{G}_{k+1}^*\hat{Q}_{k+1}^*\r^{(m+1)}_{0} \in \C^{m + t_k - k}$.
\end{enumerate}
\end{remark}

\begin{remark}
    We obtain $\vec{c}_{k + 1} = \hat{G}_{k + 1}\hat{E}(I - D_QD_Q^*)\hat{E}^*\hat{G}_{k+1}^* \hat{Q}_{k+1}^*\widetilde{\p}_{k+1}$. However, in reality, the matrix $\hat{G}_{k+1}$ is immediately updated from iteration to iteration; therefore, we store only $\hat{E}(I - D_QD_Q^*)\hat{E}^*\hat{G}_{k+1}^* \hat{Q}_{k+1}^*\widetilde{\p}_{k+1} \in \C^{m + t_k - k}$ and update it accordingly to changes in $\hat{G}_{k+1}$. So in every iteration $k$ we store $m$ vectors of size $m + t_k - k$. We recall that $t_k \le m + k$; therefore, $m + t_k - k \le 2 m$.
\end{remark}

\subsection{The algorithm}

The whole process is represented in Algorithm~\ref{alg:GMRESMRHS}

\begin{algorithm}[ht]
\label{alg::gmres_mrhs}
\begin{algorithmic}
    \Require $A$, $\vec{b}^{(i)}$, $\x_0^{(i)}, i = 1\dots M$
    
    \State $\r_0^{(1)} = \vec{b}^{(1)} - A\x_0^{(1)}$
    \State $Q_0 = \begin{bmatrix}\r_0^{(1)} / \norm{\r_0^{(1)}} \end{bmatrix}$
    \State $m = 1$, $K = 1$, $n_0 = 0, G_0 = [1]$, $H_0 = [~]$, $\hat{E}_0 = [1], E = [~], Z = [\vec{r}_0^{(1)}],  \v_1 = 1$
    \For {$k = 1, \dots$}
    \State $\p_k = Q_{k-1}G_{k-1}\hat{E}_{k-1}\v_k$
    \State $\vec{h}_k = Q_{k-1}^*A\p_k$
    \State $\q = A\p_k - Q_{k-1}\vec{h}_k$
    \If{$\norm{\q} > 0$}
    \State $Q_k = \begin{bmatrix}
        Q_{k-1} & \q / \norm{\q}
    \end{bmatrix}$,
    $H_k = \begin{bmatrix}
        H_{k-1} & \vec{h}_k \\
        \0^T & \norm{\q} 
    \end{bmatrix}$,
    $K = K + 1$
    \Else 
    \State $H_k = \begin{bmatrix}
        H_{k-1} & \vec{h}_k \\ 
    \end{bmatrix}$
    \EndIf
    \State $\begin{bmatrix}
        G_k &  R_k
    \end{bmatrix} = \textbf{qr} (H_k)$
    \State $\hat{E}_k = \begin{bmatrix}
        E & \e_k & \dots & \e_{K}
    \end{bmatrix}$
    \State $D_k = \textbf{orth}(\hat{E}_k^*G_k^*Q_k^*Z)$
    \If {$\norm{\begin{bmatrix}
        O_k & I_{K - k}
    \end{bmatrix} G_k^*Q_k^*\r_0^{(m)}} \le \varepsilon^{(m)} \norm{\vec{b}^{(m)}}$}

        \State $\x^{(m)} = \x_0^{(m)} + \sum\limits_{i = 1}^{k}\vec{y}_i^{(m)}Q_{i}G_{i}\hat{E}_i\v_i$
        \State $\vec{y}^{(m)} = R_{k}^{-1}\begin{bmatrix}
            I_{k} & O
        \end{bmatrix} G_{k}^*Q_{k}^*\r_0^{(m)}$

        \If {$(m == M)$}
            \State  \textbf{break}
        \EndIf
        \State $\r_0^{(m+1)} = \vec{b}^{(m+1)} - A\x_0^{(m+1)}$
        \State $\q = A\p_k - Q_{k-1}\vec{h}_k$
        \If {$\norm{\q} > 0$}
            \State $Q_k = \begin{bmatrix}
                Q_k & \q / \norm{\q}
            \end{bmatrix}$,
            $G_k = \begin{bmatrix}
                G_k & \0 \\
                \0^T & 1
            \end{bmatrix}$,
            $D_k = \begin{bmatrix}
            D_k  \\
            \0^T
            \end{bmatrix}$,
            $\hat{E}_k = \begin{bmatrix}
                \hat{E}_k & \0 \\ 
                \0^T       & 1        
            \end{bmatrix}$,
            $K = K + 1$
        \EndIf
        \State $\v_{k+1} = \hat{E}_k^*G_k^*Q_k^*\r_0^{(m+1)}$
        \State $E = \begin{bmatrix} E & \vec{e}_k \end{bmatrix}$, $Z = \begin{bmatrix} Z & \r_0^{(m+1)} \end{bmatrix}$, $m = m + 1$
    \Else
    \State $\v_{k + 1} = \e_{k}$
    \EndIf
    \State $\v_{k+1} = (I - D_kD_k^*)\v_{k+1}$
    \EndFor

\end{algorithmic}
\caption{GMRES for multiple right-hand sides}
\label{alg:GMRESMRHS}
\end{algorithm}

\section{Numerical results}

In this section, some numerical experiments are reported to compare the performance of the extended GMRES method with that of the extended GCR method and the seed variant of GMRES. The algorithm \ref{alg::gmres_mrhs}, as well as the competing algorithms, has been implemented using the C++ programming language. For matrix and vector operations, the BLAS and LAPACK libraries have been used.  

The problem of scattering a flat electromagnetic wave on a perfectly conducting surface has been chosen for numerical experiments. The different right-hand sides correspond to different angles of incidence of the wave. Calculations have been performed for the airplane model and for the perfectly conducting cylinder. A detailed statement of the problem and a scheme
for its numerical solution can be found in \cite{stavtsev2009application}. We highlight that the problem is solved by the method of boundary integral equations with the use of RWG basis functions. For the airplane model, the matrix $A$ of size 577680 has been formed along with 722 right-hand sides corresponding to angles from $0^\circ$ to $180^\circ$ with $0.5^\circ$ step and 2 types of polarization. For the perfectly conducting cylinder, the matrix $A$ of size 173920 has been formed along with 722 right-hand sides corresponding to angles from $0^\circ$ to $180^\circ$ with $0.5^\circ$ step.

If $\varepsilon$ is a relative tolerance and $\hat{\x}^{(i)}$ is a corresponding numerical solution, then

\[\gamma_i = \dfrac{ ||A\hat{\x}^{(i)} - \vec{b}^{(i)}||}{\varepsilon||\vec{b}^{(i)}||}\]
is the ratio between the error obtained and the required tolerance.

The comparison is represented in Tables \ref{tab:tol_2}, \ref{tab:cyl} where the maximum ratios $ \max\limits_{i = \overline{1,N}} \gamma_i $ and the geometric mean ratios $\sqrt[N]{\gamma_1 \dots  \gamma_N}$ are compared. Along with it, the total number of iterations and execution time are provided. The time of the matrix-vector product is measured apart to have an idea of the time spent on methods' operations and on matrix by vector products.

The tests have been performed on the INM RAS cluster (\url{https://cluster2.inm.ras.ru}). The single precision has been used.
\begin{table}[ht]
    \centering
    \begin{tabular}{|c|c|c|c|c|c|}
        \hline 
         method &\# iterations  & geom mean  & max rel  & time  & matvec time \\  \hline \hline
         \multicolumn{6}{|c|}{Tolerance $10^{-2}$} \\
         \hline
         GMRES & 5355 & 0.931464 & 0.999991 & 16641 & 11490 \\ 
          GCR & 5526 & 1.05426 & 1.25867 & 18323 & 12169 \\ 
         Seed & 27384 & 0.940904 & 1.00001  & 64740 & 63630 \\
         \hline

          \multicolumn{6}{|c|}{Tolerance $10^{-3}$} \\ \hline 
           GMRES & 8834 & 0.930439 & 1.00085 & 34098 & 20471 \\
           GCR & \multicolumn{5}{|c|}{not converged} \\ 
           Seed & 62215 & 0.936516 &  1.01079 & 152729 & 144565 \\
         \hline
         \multicolumn{6}{|c|}{Tolerance $10^{-4}$} \\ \hline  
         GMRES & 13654 & 1.03574 & 1.66923  & 62507  & 31727 \\ 
          GCR & \multicolumn{5}{|c|}{not converged} \\
         Seed & 115629 & 1.32348 & 1.87281  & 290379 & 268680 \\
         \hline
          
    \end{tabular}
    \caption{Performance of multi-rhs methods for the airplane model}
    \label{tab:tol_2}
\end{table}

\begin{table}[ht]
    \centering
    \begin{tabular}{|c|c|c|c|c|c|}
        \hline 
         method &\# iterations  & geom mean  & max rel  & time  & matvec time  \\  \hline \hline
         \multicolumn{6}{|c|}{Tolerance $10^{-2}$} \\
         \hline
         GMRES & 3051 &  0.895028 &  0.999979 & 1421 & 907 \\ 
         GCR & 3078 & 0.897977 & 1.00013 &  1343 &  960 \\ 
         Seed &  7900 & 0.927453  &  0.999982  & 2438 & 2315 \\
         \hline

          \multicolumn{6}{|c|}{Tolerance $10^{-3}$} \\ \hline 
          GMRES & 4111 &  0.862985 &  0.999625  & 2116 & 1169 \\
          GCR & 4116 &  1.07512 & 1.60924  & 2032 &  1217 \\ 
          Seed & 21904 & 0.856931 & 0.999939 & 7034 & 6418 \\
         \hline
         \multicolumn{6}{|c|}{Tolerance $10^{-4}$} \\ \hline 
         GMRES & 5199  &  0.874624  & 1.00114  & 2692 & 1509 \\ 
          GCR & 5181 & 51.1946 & 165.614 & 2670 &  1516\\ 
         Seed & 33264 & 0.859336 &  1.0025 &  11998 & 9746 \\
         \hline
          
    \end{tabular}
    \caption{Performance of multi-rhs methods for the perfectly conducting cylinder}
    \label{tab:cyl}
\end{table}

From the numerical results, our method is concluded to be more robust and allows us to achieve higher precision than the GCR method.

\section{Conclusion}
The Generalized Minimal Residual method has been extended for systems with the same nonsingular matrix and multiple right-hand sides that are not available simultaneously. The search space does not match the Krylov subspace anymore, since it is not rearranged after changing a right-hand side. However, the new method preserves such properties of the classical GMRES as follows: (1) a common basis for the search space and its image are constructed, (2) bases for both spaces are effectively represented within the common basis, and (3) both bases are maintained orthonormal.

The new method is mathematically equivalent to the extended GCR method \cite{lingen1999generalised}, but it is more robust, less computationally complex, and requires less memory.

The mechanism of an extension of the search space with an arbitrary residual can be applied to GMRES-like methods with flexible direction vectors.

\printbibliography
\end{document}